\documentclass{amsart}
\usepackage{graphicx}
\usepackage{latexsym}
\usepackage{hyperref}
\usepackage{amsfonts}
\usepackage[all]{xy}
\usepackage{amssymb, mathrsfs, amsfonts, amsmath}
\usepackage{amsbsy}
\usepackage{amsfonts}
\newtheorem{theorem}{Theorem}

\numberwithin{equation}{section}

\begin{document}
\title [Einstein Hermitian metrics]{Einstein Hermitian Metrics of Non Negative Sectional Curvature}
\author{Ezio de Araujo Costa }
\maketitle
\begin{center}{\bf Abstract}
\end{center}
In this paper we will prove that if $M$ is a compact simply connected Hermitian Einstein 4-manifold with non negative sectional curvature then $M$ is isometric to complex projective space $\Bbb{CP}^2$ with the Fubini-Study metric or $M$ is isometric to a product of two two-spheres $\mathbb{S}^2 \times \mathbb{S}^2$, with theirs canonical metrics.

\section{\bf Introduction}

Let $M = M^4$ be a 4-manifold. A Riemannian metric $g$ on $M$ is called {\it Einstein} if $M$ has constant Ricci curvature and called {\it Hermitian} if $g(J\cdot,J\cdot) = g(\cdot,\cdot)$ for a complex structure $J$ on $M$.
 In [5], C. LeBrun proved the following

\begin{theorem}
({\bf LeBrun)}-\label{thm1}
Let $\big( M=M^4,\,g,\, J)$ be a compact connected complex surface with metric $g$ and complex structure $J$. If $g$ is Einstein and Hermitian with respect to $J$ then only one of the following holds:
\begin{enumerate}
 \item $g$ is Kaehler-Einstein with positive Ricci curvature.
\item $M$ is isometric to $\Bbb {CP}^2 \sharp \overline{\Bbb{CP}^2}$ and $g$ is the Page metric.
\item $M$ is isometric to $\Bbb {CP}^2 \sharp \overline{\Bbb{CP}^2}$ and $g$ is the Chen-LeBrun-Weber metric.
\end{enumerate}
\end{theorem}

Using the previous theorem, C. Koca proved in [4] the following:

\begin{theorem}
\label{thm2}
({\bf Koca})-Let $\big(M=M^4,\,g,\, J)$ be a compact complex surface with metric $g$ and complex structure $J$. If $g$ is Einstein and Hermitian with respect to $J$ and $g$ has positive sectional curvature then $M$ is isometric to complex projective space $\Bbb{CP}^2$  with the Fubini-Study metric metric.
\end{theorem}

Now, consider $M$ a compact simply connected Kaehler-Einstein 4-manifold with non negative sectional curvature. In this case,  M. Berger proved in [1], that $M$ is isometric to complex projective space $\Bbb{CP}^2$ with the Fubini-Study metric or isometric to a product of two spheres $\Bbb{S}^2\times \Bbb{S}^2$, with theirs canonical metrics.

In the next two sections we will prove that the Page metric and the Chen-LeBrun-Weber metrics no has non negative sectional curvature. This will conclude the proof our main result:

\begin{theorem}
\label{thm3}
Let $\big(M=M^4,\,g,\, J)$ be a compact simply connected complex surface with metric $g$ and complex structure $J$. If $g$ is Einstein and Hermitian with respect to $J$ and $g$ has non negative sectional curvature then $M$ is isometric to complex projective space $\Bbb{CP}^2$ with the Fubini-Study metric
 or $M$ is isometric to a product of two spheres $\Bbb{S}^2\times \Bbb{S}^2$, with theirs canonical metrics.
\end{theorem}
\newpage
\section{\bf Page metric}
The Page metric (see [6]) lives in connected sum $\Bbb{CP}^2 \overline{\sharp \Bbb{CP}^2}$, where $\overline{\sharp \Bbb{CP}^2}$ is the complex projective space $\Bbb{CP}^2$ with opposite orientation.
In [4], C. Koca showed that the Page metric no has non negative sectional curvature using a computer program like {\it Maple}. In this section we will prove this result using a different argument. For this consider the page metric $g$ as in Koca [4] :

$$ g = W^2(x)dx^2 + g^2(x)(\sigma_{1}^2 + \sigma_{2}^2) + \frac{D^2}{W(x)}\sigma_{3}^2$$

where $x\in (-1,1)$,
$$W(x) = \sqrt \frac{1- a^2x^2}{[3 - a^2 - a^2(1 + a^2)x^2](1 - x^2)},$$

$$ g(x) = \frac{2}{\sqrt{3 + 6a^2 - a^4}}\sqrt{1 - a^2 x^2},$$
 $ D = \frac{2}{3 + a^2}$ and $a$ is the unique positive root of equation $ f(x) = x^4 + 4x^3 - 6x^2 + 12x - 3= 0$. Notice that  $a < 1$, since that $f(0) =-3$ and $f(1)= 8$.
\\
\\
In accord with Koca, there exists a two-plane where the sectional curvature satisfies $$K_{01} = 2[\frac{ g'W' - g''W}{gW^3}].$$

Then we have $$K_{01} = -\frac{2}{gW}F',$$ where $$F = \frac{g'}{W}.$$
\\
\\
 {\bf Claim}: There exist $c\in [0, 1)$ such that $F'(c) > 0$.
\\
\\
Proof of Claim :
\\
\\
Notice that $g' = -Ax(1 - a^2 x^2)^{-1/2}$, where $A = \frac{2a}{\sqrt{ 3 + 6a^2 - a^4}} > 0$. Moreover,
$$F =-Ax\sqrt{[3 - a^2 - a^2(1 + a^2)x^2](1 - x^2)}(1 - a^2 x^2)^{-1},$$ where $1- a^2x^2 > 0$.
Assumes that $F'(x) \leq 0$ for all $x\in [0,1)$.
\\
\\
Then $F$ is a increasing function in [0,1) and follows of this that $F(y) \leq F(x) \leq F(0) = 0$, for all $ y > x \in [0,1)$. Since that $F$ is continuo in [0,1], we have that $ 0 = F(1) \leq F(x)\leq 0$, for all $x\in [0,1)$. So $F = 0$ in [0,1) (contradiction)
\\
\\
This proves that there exists points where $K_{01} < 0$.
\newpage
\section{\bf Chen-LeBrun-Weber metric}

In [2], Chen, Lebrun and Weber proved that $ M = \Bbb{CP}^2 \overline{\sharp 2\Bbb{CP}^2}$ admits an Hermitian non-Kaehler Einstein metric $g$. In particular, there exists a Kaehler metric $h$ on $M$ of positive scalar curvature $s$ such that $g = s^{2}h$. Now consider $W_{h}^+$ and $W_{g}^+$ the self-dual Weyl part of the Weyl tensor $W$ of the respective metrics $h$ and $g$. Since that $h$ and $g$ are conformally related we have $W_{g}^+ = \frac{1}{s^2}W_{h}^+$. On the other hand, the self-dual tensor $W_{h}^+$ of the Kaehler metric $h$ has exactly two  different eigenvalues and so $W_{g}^+$ has also two different eigenvalues. By Proposition 4 of Derdzinski [3 ], $M$ admits an non trivial Killing vector field with respect to metric $g$. Assumes that $M$ has non negative sectional curvature with respect to metric $g$. Then $M$ is a compact simply connected 4-manifold with non negative sectional curvature and with a non trivial Killing vector field. By Theorem 1 of Searle and Yang in [7], we have that the Euler characteristic of $M$ satisfies $\chi(M) \leq 4$ which contradicts the fact of that $\chi(\Bbb{CP}^2 \overline{\sharp 2\Bbb{CP}^2}) = 5$.

Author's address:

Mathematics Department,  Federal University of Bahia,

zipcode: 40170110- Salvador -Bahia-Brazil

Author's email

 ezio@ufba.br
\end{document}